\documentclass[11pt, final]{amsproc}
\usepackage[cp850]{inputenc}
\usepackage{amssymb}
\usepackage[mathscr]{eucal}

\usepackage[all]{xy}
\usepackage{amscd}
\usepackage{amsmath}
\usepackage{setspace}
\newcommand{\R}{\mathbb R}

\def\dim{\operatorname{dim}}
\def\supp{\operatorname{supp}}
\newcommand{\Z}{\mathscr{X}}
\newcommand{\KP}{{\sf K}\hspace{-1pt}{\sf P}}

\theoremstyle{plain}
\newtheorem{theorem}{Theorem}
\newtheorem{proposition}{Proposition}
\newtheorem{lemma}{Lemma}
\newtheorem{corollary}{Corollary}

\theoremstyle{remark}

\title{The Kalton-Peck space is the complexification \\ of the real Kalton-Peck space}

\author[J.M.F. Castillo]{Jes\'us M. F. Castillo}
\address{Universidad de Extremadura, Instituto de Matem\'aticas Imuex,
E-06011 Badajoz, Spain.}
\email{castillo@unex.es}

\author{Yolanda Moreno}
\address{Universidad de Extremadura, Instituto de Matem\'aticas Imuex,
E-07001 Cßceres, Spain.}
\email{ymoreno@unex.es}

\thanks{This research has been supported in part by project PID2019-103961GB funded by MINCIN and project IB20038 funded by Junta de Extremadura.}

\begin{document}

\begin{abstract} The Kalton-Peck $Z_2$ space is the derived space obtained from the scale of $\ell_p$ spaces by complex interpolation at $1/2$.
If we denote by  by $Z_2^{real}$ the derived space obtained from the scale of $\ell_p$ spaces by real interpolation at $(1/2, 1/2)$, we show that $Z_2$ is the complexification of $Z_2^{real}$. We also show that $Z_2^{real}$ shares the most important
properties of $Z_2$: it is isomorphic to its dual, it is singular and contains no complemented copies of $\ell_2$.\end{abstract}
\maketitle

This paper outgrowths from \cite{cck}, in which, in the words of its authors, \emph{it is shown that Rochberg's generalized interpolation spaces $\Z^{(n)}$ arising from analytic families of Banach spaces form exact sequences $0\to  \Z^{(n)} \to \Z^{(n+k)} \to \Z^{(k)} \to 0$ and that nontriviality, having strictly singular quotient map, or having strictly cosingular embedding depend only on the basic case
$n=k=1$}. Our attempt has been to bring down to sound earth its Open end 6.3: \emph{It is the feeling of the authors that most of the work done in this paper could be reproduced for real interpolation by either the $K$
or $J$ methods with a careful analysis of the work done in [Carro et al. 1995]. It
would be interesting to know to what extent the same occurs for other interpolation
methods.}\medskip

That is what we will do in this first part of the paper in the following way: It is a fact well known to all those who know it that most interpolation methods generate exact sequences of the interpolated spaces; say, if $[X_0,X_1]_\mu$ is the interpolated space obtained from the pair $(X_0, X_1)$ with parameters set at $\mu$, then there exists a natural exact sequence
$\xymatrix{0\ar[r]& [X_0,X_1]_\mu\ar[r]& \mathcal R^{(2)} \ar[r]& [X_0,X_1]_\mu \ar[r]& 0}$ in which $[X_0,X_1]_\mu$ is the first Rochberg space $\mathcal R^{(1)}$ and $\mathcal R^{(2)} $ will be the second Rochberg derived space. After that, the higher order Rochberg derived spaces
\cite{rochberg} can be rather naturally generated when the method provides a sequence of \emph{compatible} interpolators (see \cite{racsam,sym}) and form natural exact sequences $\xymatrix{0\ar[r]& \mathcal R^{(n)}  \ar[r]& \mathcal R^{(n+m)} \ar[r]& \mathcal R^{(m)} \ar[r]& 0}$ (see \cite{cck,ccc}). Once again, most standard interpolation methods provide associated sequences of compatible interpolators. This is more or less implicit in the papers of Cwikel et al. \cite{cjmr}, Carro et al. \cite{caceso} and Rochberg \cite{rochberg} and made explicit when Cwikel, Kalton, Milman and Rochberg introduce their unifying method \cite{CKMR}, from now on called the CKMR method.

On the other hand, the theory created by Kalton \cite{kalt} establishes that exact sequences of quasi Banach spaces are in correspondence with
a special type of nonlinear maps called quasi linear maps. Given an exact sequence
$\xymatrix{0\ar[r]& X \ar[r]& X_\Omega \ar[r]& X \ar[r]& 0}$ with associated quasilinear map $\Omega$ and another
$\xymatrix{0\ar[r]& X \ar[r]& X_\Phi \ar[r]& X \ar[r]& 0}$ with associated quasilinear map $\Phi$, the two exact sequences are called projectively equivalent \cite{kaltpeck,hmbst}
if there is an scalar $\lambda$ such that the diagram
$$\xymatrix{0\ar[r]& X\ar[d]_\lambda \ar[r]& X_\Omega \ar[r]& X \ar[r]\ar@{=}[d]& 0\\
0\ar[r]& X \ar[r]& X_\Phi \ar[r]& X \ar[r]& 0}$$ is commutative, which means that there is a linear map $L$ such that $\lambda \Omega - \Phi$ is the sum of a bounded and a linear map, both $X\to X$.

The same occurs for the exact sequences of Rochberg spaces. In his case, these quasilinear maps have special form and additional properties (see below) \cite{kaltpeck,racsam} and are called differentials \cite{ccfg,jussieu}, which justifies that we call the process of obtaining the associated differential out from an interpolation method as derivation. The derivation process for the CKRM method is described and studied in detail in \cite{CKMR}. We make a rapid survey of what the reader needs to know about the CKRM method to ease the reading of this paper.

\section{Preliminaries:. The CKRM method}

Let $\mathbf{Ban}$ be the class of all complex Banach spaces. A mapping $\mathcal X:\mathbf{Ban} \to \mathbf{Ban}$ will be called a pseudolattice if
\begin{enumerate}
\item[(i)] for each $B\in \mathbf{Ban}$ the space $\mathcal X(B)$ consists of $B$-valued sequences $\{b_n\}_{n\in \mathbb Z}$ and if
\item[(ii)] whenever $A$ is a closed subspace of $B$ it follows that $\mathcal X(A)$ is a closed subspace of $\mathcal X(B)$ and if
\item[(iii)] there exists a positive constant $C>0$ such that, for all $A,B\in \mathbf{Ban}$, and all bounded linear operators $T:A\to B$ and every sequence $\{a_n\}_{n\in \mathbb Z}\in \mathcal X(A)$, the sequence $\{T(a_n)\}_{n\in \mathbb Z}\in \mathcal X(B)$ satisfies the estimate $$\| \{T(a_n)\}_{n\in \mathbb Z} \|_{\mathcal X(B)}\leq C(\mathcal X)\|T\|_{A\to B}\|\{a_n\}_{n\in \mathbb Z}\|_{\mathcal X(A)} .$$
\end{enumerate}

Fix a pair of pseudolattices  $\mathbf {X}=\{ \mathcal X_0, \mathcal X_1\}$. Given a compatible pair of Banach spaces $B=(B_0,B_1, \Sigma)$ ---i.e., a pair of spaces $B_0, B_1$ considered as subspaces of another Banach space $\Sigma$--- we define $\mathcal J(\mathbf {X},B)$ to be the space of all $(B_0\cap B_1)$-valued sequences $\{b_n\}_{n\in \mathbb Z}$ for which the sequence $\{e^{jn}b_n\}_{n\in \mathbb Z}\in \mathcal X_j(B_j) $ for $j=0,1$. This space is normed by
$$\| \{b_n\}_{n\in \mathbb Z}\|_{\mathcal {J} (\mathbf {X},B)}=\max_{j=0,1} \| \{e^{jn}b_n\}\|_{\mathcal X_j(B_j)}.$$

The pseudolattice pair $\mathbf {X}$ is \textit{nontrivial} if, for the special one-dimensional Banach pair $(\mathbb C, \mathbb C)$ and each $z\in \mathbb A=\{z\in \mathbb C:1<|z|<e \}$ (the open annulus) there exists $\{b_n\}_{n\in \mathbb Z}\in \mathcal {J} (\mathbf {X},(\mathbb C, \mathbb C)) $ such that the series $\sum_{n\in \mathbb Z}z^nb_n$ converges to a nonzero number. The pseudolattice pair $\mathbf {X}$ is \textit{Laureant compatible} if it is nontrivial and for every $z\in \mathbb A$ the Laureant series $\sum_{n\in \mathbb Z} z^nb_n$ converges absolutely with respect to the norm of $B_0+B_1$. Therefore the sum of this series is an analytic function of $z$ in $\mathbb A$ and can be differentiated term-by-term. The series for its derivative $\sum_{n\in \mathbb Z} nz^{n-1}b_n$ also converges absolutely in $B_0+B_1$. See both claims at \cite[page 248]{CKMR}.\\

Given a compatible pair $B=(B_0,B_1)$ of Banach spaces and $0<\theta<1$, we define the \textit{interpolation space} $B_{\mathbf {X}, \theta}$ to consist of all elements of the form $b=\sum_{n\in \mathbb Z}e^{\theta n}b_n$ with $\{b_n\}_{n\in \mathbb Z}\in \mathcal J(\mathbf {X},B)$, endowed with the natural quotient norm:
$$\| b \|_{B_{\mathbf {X}, \theta}}=\inf \left\lbrace \| \{b_n\}_{n\in \mathbb Z}\|_{\mathcal J(\mathbf {X},B)}: b= \sum_{n\in \mathbb Z}e^{\theta n}b_n   \right\rbrace .$$
According to the previous claims one may ``think" every element $\{b_n\}_{n\in \mathbb Z}\in \mathcal J(\mathbf {X},B)$ as the analytic map $\sum_{n\in \mathbb Z} z^nb_n$, where $z\in \mathbb A$, with all the precautions. This is, we shall informally write $\{b_n\}_{n\in \mathbb Z}=\sum_{n\in \mathbb Z} z^nb_n.$  Therefore, we have the natural \textit{evaluation map} $\delta_{\theta}:\mathcal J(\mathbf {X},B) \longrightarrow B_{\mathbf {X}, \theta}$ given by the rule $\delta_{\theta}(\{b_n\}_{n\in \mathbb Z})=\sum_{n\in \mathbb Z} e^{\theta n} b_n.$\\

\subsection*{Intermission: The module structure} Banach spaces $Z$ with a $1$-unconditional basis admit an obvious structure of $\ell_\infty$-module given by the product $\ell_\infty \times Z\to Z$ in which $(\xi x)(n)= \xi(n)x(n)$. We are especially interested in the situation in which one deals with pairs $B=(B_0,B_1)$ of spaces with a joint $1$-unconditional basis. In such case, we say a pseudolattice pair $\mathbf X$ admits an $\ell_{\infty}$-module structure if there is $C>0$ such that for every $\{b_n\}_{n\in \mathbb Z}\in \mathcal J(\mathbf {X},B)$ and $a\in \ell_{\infty}$, the following holds:
\begin{enumerate}
\item $\{a\cdot b_n\}_{n\in \mathbb Z}\in \mathcal J(\mathbf {X},B)$
\item $\|\{a\cdot b_n\}_{n\in \mathbb Z}\|_{\mathcal J(\mathbf {X},B)}\leq C\|a\|_{\infty}\|\{b_n\}_{n\in \mathbb Z}\|_{\mathcal J(\mathbf {X},B)}$
\end{enumerate}
The natural example is $\mathbf X=\{\ell_{p_0}, \ell_{p_1}\}$ because one trivially has $\|a\cdot b_n\|_{B_j}\leq \|a\|_{\infty} \|b_n\|_{B_j},\;\;j=0,1$; and thus, $\|\{a\cdot e^{jn} b_n\}_{n\in \mathbb Z}\|_{\ell_{p_j}(B_j)} \leq \|a\|_{\infty}\| \{e^{jn} b_n\}_{n\in \mathbb Z}\|_{\ell_{p_j}(B_j)}$ for $j=0,1$.\textbf{End of the intermission}\\

Given  $C \geq 1$ a $C$-extremal for a given $b\in B_{\mathbf {X}, \theta}$ is a sequence $\{b_n\}_{n\in \mathbb Z}$ so that $\delta_{\theta}(\{b_n\}_{n\in \mathbb Z})=b$ and $\|\{b_n\}_{n\in \mathbb Z}\|\leq C\|b\|$. We write $S(b)=\{b_n\}_{n\in \mathbb Z}$ and say that $S$ is $C$-\textit{bounded selector} for the map $\delta_{\theta}$. The associated \textit{differential map} $\Omega$ is defined by the $\Sigma$-valued map
\begin{equation}\label{derivat}
 \Omega(b) =\delta_{\theta}'S(b)=\sum_{n\in \mathbb Z} ne^{\theta(n-1)}b_n.
\end{equation}
Observe that $\Omega$ does not take (necessarily) values in $B_{\mathbf {X}, \theta}$. We will work under the condition that $\mathbf {X}$ \textit{admits differentiation},  \cite[Lemma 3.11]{CKMR} which is a rather technical condition. However, the reader may keep in mind that $\mathbf {X}$ admits differentiation if the shift operator is an isometry on $\mathcal X_j(B_j)$ for $j=0,1$, see \cite[Lemma 3.6.]{CKMR}. This will be our case.\medskip

If $\mathbf X$ is a pseudolattice pair with $\ell_{\infty}$-module structure and $B=(B_0,B_1)$ is a pair of spaces with a joint $1$-unconditional basis, the corresponding differential map $\Omega$ for $B_{\mathbf X, \theta}$ is a \textit{centralizer} in the sense of Kalton \cite{kaltdiff}, which means that $\Omega(\xi x) - \xi \Omega (x) \in B_{\mathbf {X}, \theta}$ for all $\xi\in \ell_\infty$ and $x\in B_{\mathbf {X}, \theta}$
and $\|\Omega(\xi x) - \xi \Omega (x)\|\leq Q \|xi\|_\infty\|x\|$. This follows easily from the fact that the pseudolattice $\{\ell_{p_0},\ell_{p_1}\}$ has an $\ell_{\infty}$-module structure.\\

Let as before $B$ be a compatible couple with ambient space $\Sigma$. The derived space $dB_{\mathbf X,\theta}$ is the set of couples $(x,y)\in \Sigma \times B_{\mathbf X,\theta}$ for which the following quasi-norm  $$\|x-\Omega(y)\|_{B_{\mathbf X,\theta}}+\|y\|_{B_{\mathbf X,\theta}}$$
makes sense and is finite, where $\Omega$ is the differential (\ref{derivat}). If $B_{\mathbf X,\theta}$ contains no (uniform) copies of $\ell_1^n$ for every $n\in \mathbb N$ then such a quasi-norm is equivalent to a norm by a result of Kalton \cite{kalt}. In our case, $B_{\mathbf X,\theta}=\ell_2$ so the condition is trivially satisfied.\\

\subsection*{The CKMR method and the real interpolation method.}

Even the meaning of ``the real interpolation method" is somewhat ambiguous since there are many \emph{real} methods. This is not a problem since one usually uses the $K-$ and $J-$ methods, and these are equivalent \cite[Chapter 3, 3.3]{berglof}. In this paper however we will need to rely on the original ``espaces de moyennes" real method of Lions and Peetre \cite{LiPe} out of which the $K-$ and $J-$ methods were derived. This method involves a discrete and a continuous version and four parameteres ($\xi_0,\xi_1, p_0, p_1$), that we will set at
$p_0=1, p_1=\infty, \xi_0=1, \xi_1=-1$ for the pair $(\ell_1, \ell_\infty$) and that when adequately fixed produce and equivalent method to the standard real interpolation \cite[Chapter 3,3.12]{berglof}. Regarding the equivalence with the CKMR method, the pseudolattice corresponding to the real method correspond to the choice $\mathcal X_j=\ell_{p_j}$ and admits differentation since the shift operator is clearly an isometry on $\ell_{p_j}(B_j)$ for $j=0,1$, see \cite[Lemma 3.6.]{CKMR}.  The equivalence with the CKRM method in practice  means to multiply by the weight $e^{-n\theta}$. See the discussion in \cite[Paragraphs 2 and 4, page 251]{CKMR}.\\

The advantage of using the CKRM method is the explicit existence of differentials with the manageable simple formula (\ref{derivat}), something that is much harder to obtain from the descriptions in either \cite{cjmr} or \cite{caceso}. The real interpolation method applied to the pair $(\ell_\infty, \ell_1)$ thus produces the interpolation space $(\ell_\infty, \ell_1)_{1/2, 1/2}=\ell_2$ with derived twisted Hilbert space $Z_2^{real}$, that will be referred to as \emph{the real} Kalton-Peck space.

\section{The real Kalton-Peck space}\label{kaltonpeque}

Complex interpolation applied to the pair $(\ell_\infty, \ell_1)$ produces $\ell_2$ as interpolated space and the celebrated Kalton-Peck space $Z_2$ as derived space. The real method applied to the pair $(\ell_\infty, \ell_1)$ produces, as we have just said, $\ell_2$ as interpolated space. To determine the derived space $Z_2^{real}$ what we will do is to approach the real method as a CKMR method.

\begin{proposition}(Lions-Peetre) The differential associated to $(\ell_{p_0}, \ell_{p_1})_{\theta, p}=\ell_p$ is \begin{equation}\label{kpreal}
\KP^{real}(a)=e^{-\theta}\sum_m - \left(\frac{p}{p_0}-\frac{p}{p_1}\right) \left[ \log \frac{|a_m|}{\|a\|} \right] a_me_m ,
\end{equation}
for $a\in \ell_p$ and $\frac{1}{p}=\frac{1-\theta}{p_0}+\frac{\theta}{p_1}$, $0<\theta<1$ and $1\leq p_0,p_1<\infty$. Here $[ \cdot]$ means ``the entire part of".
\end{proposition}
\begin{proof}
That $\ell_p=(\ell_{p_0},\ell_{p_1})_{\theta,p}$ is established in \cite[Theorem (I.I), Chapitre VII]{LiPe}.
With the same language and notation of the paper, observe that the starting point for the proof of this is that the space of moyennes is $S(p_0, \xi_0, \R; p_1, \xi_1, \R)= \R$ with the decomposition
$a= (\dots, 0, a, 0, \dots)$, namely $a(m)= \sum_{-\infty}^{+\infty} w_n(m)$ with $w_0(m)=a(m)$ and $w_n(m)=0$ otherwise.

Let us work with the CKMR method, fix $C>1$ and obtain a $C$-extremal; namely, given $a=(a(m))_m$ we look for $\{b_n\}_{n\in \mathbb Z}$ such that \begin{itemize}
\item[(a)] $\delta_{\theta}(\{b_n\}_{n\in \mathbb Z})= a$.
\item[(b)] $\|\{b_n\}_{n\in \mathbb Z}\|\leq C \|a\|_{\ell_p}.$\end{itemize}
The proof of \cite[Theorem (I.I), Chapitre VII]{LiPe} contains the idea to obtain the $C$-extremals: set $\lambda=\frac{p}{p_0}-\frac{p}{p_1}$ so that
\begin{equation*}
\left\{ \begin{array}{lcc}
            p_0(1+\lambda \theta)=p, \\
             p_1(1-\lambda(1-\theta))=p,\\
             \end{array}
   \right.
\end{equation*}
$\|a\|=1$.  Define $\{b_n\}_{n\in \mathbb Z}$ by $b_n(m) = w_{n + [\lambda \log |a(m)|]}(m)$, which yields
$b_n(m) = w_{0}(m) = a(m)$ when $n= - [\lambda \log |a(m)|]$ and $0$ otherwise. After the corresponding multiplication by the weight $e^{-n\theta}$ above mentioned we get:
\begin{equation}\label{b}
b_n(m)= \left\{ \begin{array}{lcc}
             e^{-n\theta}a_m &,& n=-[\lambda \log |a_m|] \\
             \\ 0 &,& \mathrm{otherwise.}\\
             \end{array}
   \right.
\end{equation}
It is clear  that (a) holds. We check (b):
\begin{eqnarray*}
\|\{b_n\}_{n\in \mathbb Z}\|_{\ell_{p_0}(\ell_{p_0})}^{p_0}&=&\sum_{n} \|b_n\|^{p_0}_{\ell_{p_0}}\\
&=& \sum_{n=-[\lambda \log |a_m|]} \left |e^{-n\theta}a_m \right|^{p_0}\\
& \leq & \sum_m e^{-\theta(-\lambda \log |a_m|)p_0} |a_m|^{p_0}\\
&=& \sum_m  |a_m|^{p_0+\lambda \theta p_0}\\
&=&\sum_m |a_m|^p=1.
\end{eqnarray*}
Analogously, $\|\{e^nb_n\}_{n\in \mathbb Z}\|_{\ell_{p_1}(\ell_{p_1})}\leq 1$ and thus this element is a $C$-extremal. Now, the differential corresponding to the CKMR method will be $\KP^{real}(a)=\sum_{n\in \mathbb Z} n e^{\theta(n-1)}b_n$, namely
\begin{equation}\label{kpreal}
\KP^{real}(a)=e^{-\theta}\sum_m - \left(\frac{p}{p_0}-\frac{p}{p_1}\right) \left[ \log \frac{|a_m|}{\|a\|} \right] a_me_m.\qedhere
\end{equation}
\end{proof}

We prove now the result in the title:

\begin{theorem} The space $Z_2$ is the complexification of $Z_2^{real}$. \end{theorem}
\begin{proof} Let us denote by $\ell_2(\R)$ the real infinite dimensional separable Hilbert space and by $\ell_2$ the complex Hilbert space. We show first that there is a commutative diagram

$$\xymatrix{
0\ar[r]& \ell_2 \ar[r] & Z_2\ar[r]&\ell_2 \ar[r]&0\\
0\ar[r]& \ell_2(\R) \ar[u]^{e^{\theta}\;\imath}\ar[r]& Z_2^{real}\ar@{-->}[u]\ar[r] &\ell_2(\R) \ar[r]\ar[u]_\imath&0}$$
in which the vertical arrow $\imath$ is plain inclusion of the real $\ell_2(\R)$ into the complex $\ell_2$. Indeed
\begin{eqnarray*}
\left\| \left(\KP \imath - e^{\theta}\imath \KP^{real}\right)a\right\|_2 &=& \left\| 2\sum_m \log \frac{|a_m|}{\|a\|} a_me_m  - 2e^\theta e^{-\theta}\sum_m \left[ \log \frac{|a_m|}{\|a\|} \right] a_me_m\right\|_2\\
&=& \left\| 2\sum_m \left( \log \frac{|a_m|}{\|a\|} - \left[ \log \frac{|a_m|}{\|a\|} \right] \right) a_me_m\right\|_2\\
&\leq & 2 \left\| \sum_m \frac{1}{2} a_me_m \right\|_2 =  \|a\|_2.\end{eqnarray*}
We have thus shown that $\KP \imath - e^{\theta}\imath \KP^{real}: \ell_2(\R)\longrightarrow \ell_2(\R)$ is a bounded map, from which assertion about the diagram follows. The vertical middle dotted arrow exists by the nature of the diagram.\\

We now invoke the results in \cite[Section 4]{complex} where it is shown that the complexification of the real twisted Hilbert space
generated by the map $\KP^r(x) = x\log \frac{\|x\|}{|x|}$ on $\ell_2(\R)$ is the Kalton-Peck $Z_2$ space
generated by $\KP$ (and the same occurs after multiplying $\KP$ by any nonzero scalar).  We have thus shown that $\KP|_{\ell_2(\R))} = \KP^r$ and $\KP^{real}$ are projectively equivalent. Thus, $Z_2$ is a complexification of $Z_2^{real}$. \end{proof}

We now translate three important properties from $Z_2$ to $Z_2^{real}$ while adding a few original features to the proofs, something that could be interesting on its own. Recall that an exact sequence $\xymatrix{0\ar[r]& X \ar[r]& X_\Omega \ar[r]& X \ar[r]& 0}$ with associated quasilinear map $\Omega$ (or the quasilinear map $\Omega$) is called \emph{singular} when the quotient map is a strictly singular operator. This occurs \cite{strict} if and only if the restriction of $\Omega$ to any infinite dimensional subspace is never the sum of a bounded plus a linear map \cite{hmbst}.

\begin{proposition} $\;$\begin{itemize}
\item[\rm(1)] $Z_2^{real}$ is isomorphic to its dual.
\item[\rm(2)] $\KP^{real}$ is singular.
\item[\rm(3)] $Z_2^{real}$ does not contain complemented copies of $\ell_2$
\end{itemize}
\end{proposition}
\begin{proof} Assertion (1) is consequence of the Kalton-Peck inequality $|xy\log\frac{x}{y}|\leq C xy$ for positive $x,y>0$ (see \cite{kaltpeck}) which means that the quasilinear map that defines the dual sequence
$\xymatrix{0\ar[r]& \ell_2^*\ar[r]& Z_2^* \ar[r]& \ell_2^*\ar[r]& 0}$ is $-\KP$. The same occurs to $\KP^{r}$, hence to $\KP^{real}$.
The details, that we sketch next, can be found with different levels of depth in \cite{kaltpeck,ccc,sym,dual}. We follow \cite{dual}: Two quasiinear maps $\Omega: B \to A$ and $\Phi: A^*\to B^*$ are called bounded duals one of the other if there is $C>0$ such that for every $b\in B, a^*\in A^*$ one has $$ |\langle \Omega b, a^*\rangle + \langle b, \Phi a^*\rangle |\leq C \|b\|\|a^*\|.$$
If $\Omega$ generates the exact sequence $0\to A \stackrel{\imath}\to A\oplus_\Omega B \stackrel{\pi}\to B \to 0$ and $\Phi$ is a bounded dual of $\Omega$ then $\Phi$ generates the dual sequence $0\to B^*\to (A\oplus_{\Omega} B)^* \to A^* \to 0$, with the meaning that there is an operator
$D: B^*\oplus_\Phi A^* \longrightarrow (A\oplus_\Omega B)^*$ given by$$\langle D(b^*, a^*), (a,b)\rangle = \langle b^*, b\rangle + \langle a^*, a \rangle$$
making a commutative diagram
$$\xymatrix{
0\ar[r] &B^*\ar[r]^-{\imath}\ar@{=}[d] & B^*\oplus_\Phi A^* \ar[r]^-{\pi}\ar[d]^D & A^*\ar[r]\ar@{=}[d] &0\\
0\ar[r] &B^*\ar[r]_-{\pi^*} & (A\oplus_\Omega B)^* \ar[r]_-{\imath^*} & A^*\ar[r]&0}$$
Finally, the Kalton-Peck inequality above shows that $-\KP$ is a bounded dual of $\KP$.\medskip

To prove (2) it is enough to show that there no subspace $W\subset \ell_2$ so that the restriction $\KP^{real}|_W$ is bounded. And this occurs, by the transfer principle \cite{strict,ccs,hmbst} if and only if there is no sequence $(u_n)$ of consecutive blocks so that
the restriction of $\KP^{real}$ to the subspace $W=[(u_n)_n]$ spanned by those blocks is bounded. We work with $\KP$, which is slightly simpler, with an argument which is immediately valid for $\KP^r$, hence for $\KP^{real}$.\medskip

Normalize the elements $u_n$ in $c_0$ so that $\|u_n\|_\infty=1$. Let $|x|$ denote the size of a finite element $x$, i.e. the cardinal of its support, and set $u_n = \sum_{j\in F_n} \lambda_{n,j} e_j$ with $|F_n|=|u_n|$. The element $\sum^N u_n$ is such that $\|\sum^N u_n\|_\infty=1$, $\|\sum^N u_n\|_1= \sum^N \|u_n\|_1$ and $\|\sum^N u_n\|_2 = \left(\sum^N \|u_n\|_2^2\right)^{1/2}$. The holomorphic function
$$f(z) = \left(\sum_{n=1}^{n=N} \sum_{j\in F_n} |\lambda_{n,j}|^{2z}\right)^{(1/2)(2z -1)}\sum^N u_n$$
is obviously a selector for $\sum^N u_n$ since $f(1/2) = \sum^N u_n$. We estimate its norm. Since
$$f(1)= \left(\sum_{n=1}^{n=N} \sum_{j\in F_n} |\lambda_{n,j}|^{2}\right)^{1/2}\sum^N u_n = \left(\sum^N \|u_n\|_2^2\right)^{1/2}\sum^N u_n$$
and thus $\|f(1)\|_\infty = \left( \sum^N \|u_n\|_2^2 \right)^{1/2} = \left\|\sum^N u_n\right\|_2$, the same occurs on points $1+it$. On the other hand, since
$$f(0) = \left(\sum_{n=1}^{n=N} \sum_{j\in F_n} 1\right)^{-1/2}\sum^N u_n = \left(\sum_{n=1}^{n=N} |u_n| \right)^{-1/2}\sum^N u_n$$
and thus
\begin{eqnarray*} \|f(0)\|_1 &=& \frac{\left \|\sum^N u_n \right\|_1}{ \left(\sum^N |u_n|\right)^{1/2}} = \frac{\sum^N \|u_n\|_1}{ \left(\sum^N |u_n|\right)^{1/2}}\\
&=& \frac{\sum^N \sum_{j\in F_n} |\lambda_{n,j}|}{ \left(\sum^N |u_n|\right)^{1/2}} \leq \frac{\sum^N |u_n|^{1/2} \|u_n\|_2}{ \left(\sum^N |u_n|\right)^{1/2}}\\
&\leq& \frac{\left( \sum^N |u_n|\right)^{1/2} \left(\sum^N \|u_n\|_2^2\right)^{1/2}}{ \left(\sum^N |u_n|\right)^{1/2}}\\
&=& \left(\sum^N \|u_n\|_2^2\right)^{1/2} =\left\|\sum^N u_n\right\|_2\end{eqnarray*}
the same occurs on points $e^{it}$. This yields that $\Omega\left( \sum^Nu_n \right) = f'(1/2)$ is an acceptable differential. Now, if
$$g(z) = \left(\sum_{n=1}^{n=N} \sum_{j\in F_n} |\lambda_{n,j}|^{2z}\right)^{(1/2)(2z -1)}= \left(\sum_{n,j } |\lambda_{n,j}|^{2z}\right)^{z - 1/2)}$$
then $\log g(z) = (z-\frac{1}{2}) \log \left( \sum_{n,j} |\lambda_{n,j}|^{2z}\right)$
and thus
$$\frac{g'(z)}{g(z)} = (z-\frac{1}{2}) \frac{\sum_{n,j} 2|\lambda_{n,j}|^{2z} \log |\lambda_{n,j}}{\sum_{n,j} |\lambda_{n,j}|^{2z}}
+ \log \left( \sum_{n,j} |\lambda_{n,j}|^{2z}\right)$$
which yields
$$\frac{g'(1/2)}{g(z)} = \log \left( \sum_{n,j} |\lambda_{n,j}|\right) = \log \sum_n \|u_n\|_1$$
Therefore
$$f'(1/2) = \log \sum_n \|u_n\|_1 \left(\sum_n u_n\right).$$

This implies that $\Omega|_W$ is unbounded. The same occurs to $\KP^r$ and $\KP^{real}$. Consequently, $\KP^{real}$ is singular.\medskip

To prove (3) we make a detour.

\begin{lemma} If $\xymatrix{0\ar[r]& \ell_2 \ar[r]^\jmath& Z\ar[r]^Q&\ell_2 \ar[r]&0}$ is an exact sequence with $Q$ strictly singular and $\jmath$ strictly cosingular then $Z$ does not contain complemented copies of $\ell_2$\end{lemma}
\begin{proof} Assume then that $Z$ contains a complemented copy $B$ of $\ell_2$ and let $P: Z\to B$ be a projection. Either
$\ker Q_{|B}$ is finite or infinite dimensional. If it is finite dimensional, $Q$ is an isomorphism on some infinite dimensional subspace of
$Z$, which is impossible. If it is infinite dimensional, let $P': B\to \ker Q_{|B}$ be a continuous linear projection, that exists because $B$ is Hilbert. Thus $P'P: Z \to \ker Q_{|B}$ is a projection, as well as $P'P_{| ker Q}$, which is once more impossible since
$\jmath$ is cosingular\end{proof}

Thus, if $\xymatrix{0\ar[r]& \ell_2 \ar[r]^\jmath& Z\ar[r]^Q&\ell_2 \ar[r]&0}$ is a singular exact sequence for which there is a commutative diagram
$$\xymatrix{
0\ar[r]& \ell_2\ar[d]_\alpha \ar[r]^{\jmath}& Z\ar[r]^Q&\ell_2 \ar[r]&0\\
0\ar[r]& \ell_2 \ar[r]_{Q^*}& Z^*\ar[r]_{\jmath^*}&\ell_2 \ar[r]\ar[u]_\gamma&0}$$
in which $\alpha, \gamma$ are isomorphisms then $Z$ cannot contain complemented copies of $\ell_2$ because an operator $T$ such that $T^*$ is strictly singular must be strictly cosingular. Thus, the sequence $\xymatrix{0\ar[r]& \ell_2(\R) \ar[r]& Z_2^{real}\ar[r]&\ell_2(\R) \ar[r]&0}$ is singular and cosingular, so the previous lemma applies. \end{proof}

\section{Concluding remarks}

Having singular differential is a rather demanding condition. For instance:

\begin{proposition}\label{prop} Let $(X, X^*)$ be a real or complex interpolation pair of Banach spaces with a common unconditional basis for which there exists a continuous inclusion $X^*\to X$ and such that $(X^*, X)_{1/2}=\ell_2$ with differential $\Omega$. If $\Omega$ is singular then $X^*$ does not contain $\ell_2$.\end{proposition}
\begin{proof} Indeed, let $(u_n)$ be blocks in $X^*$ so that $[(u_n)] \simeq \ell_2$. Pick $u\in [(u_n)]$ and observe that
\begin{eqnarray*} \|u\|_{X} &=& \sup \{ <y,u> : \|y\|\leq 1;\; y\in X^*\}\\ &\geq& \sup \{ <y,u> : \|y\|\leq 1: \; y\in [(u_n)]\}\\ &=& \|u\|_{[(u_n)]^*} \sim \|u\|_2.\end{eqnarray*}
Since $\|u\|_{X}\leq \|u\|_2$ it turns out that $\|u\|_{X^*} \sim \|u\|_2$. Thus, the norms of $X$ and $X^*$ are equivalent on $[(u_n)]$, and this obliges $\Omega|_{[(u_n)]}$ to be bounded.\end{proof}

Optimistic readers could now easily believe the following conjecture:\medskip

\emph{Conjecture. Let $(X, X^*)$ be an interpolation pair of Banach spaces with a common unconditional basis for which there exists a continuous inclusion $X^*\to X$ and such that $(X^*, X)_{1/2}=\ell_2$ with differential $\Omega$. If $\Omega$ is singular then $X$ and $X^*$ are incomparable.}\medskip

We do not have a proof for that. Optimistic readers should be warned: an example in \cite{ccs} provides two incomparable spaces $(A,B)$ not containing $\ell_2$ so that the complex differential $\Omega$ at $1/2$ is an isomorphism on a complemented copy of $\ell_2$. Thus,
the converse for the assertion above and for that in Proposition \ref{prop} fail.

\end{document}